\documentclass[10pt,a4paper,twoside,reqno]{amsart}

\usepackage{amsfonts}
\usepackage{amssymb}
\usepackage{amsmath}    %multi-line etc. Subsumes ams-text,amsbsy,amsopn
\usepackage{amsthm}     %provide proof env
\usepackage{amscd}      %commutative diagrams
\usepackage{euscript}
\usepackage[english]{babel}

\newcommand{\NN}{{\mathbb N}}
\newcommand{\ZZ}{{\mathbb Z}}

\newcommand{\QQ}{{\mathbb Q}}
\newcommand{\CC}{{\mathbb C}}

\newcommand{\PP}{{\mathbb P}}
\newcommand{\GG}{{\mathbb G}}

\newcommand{\PPn}{{\PP^n}}

\newcommand{\PV}{{\PP(V)}}

\newcommand{\Extl}{{{\mathcal E\!xt}}}

\newcommand{\OPN}{{\mathcal O}_{\PPn}}

\DeclareMathOperator{\HH}{H}

\DeclareMathOperator{\Ext}{Ext}
\DeclareMathOperator{\Hom}{Hom}
\DeclareMathOperator{\Aut}{Aut}

\DeclareMathOperator{\SL}{SL}
\DeclareMathOperator{\GL}{GL}

\DeclareMathOperator{\coker}{Coker}

\DeclareMathOperator{\Id}{Id}
\DeclareMathOperator{\Stab}{Stab}
\DeclareMathOperator{\codim}{codim}
\DeclareMathOperator{\rank}{rank}
\DeclareMathOperator{\Supp}{Supp}

\DeclareMathOperator{\hd}{hd}           %homological dimension

\newtheorem{thm}{Theorem}[section]
\newtheorem{lemma}[thm]{Lemma}
\newtheorem{prop}[thm]{Proposition}

\newtheorem{corol}[thm]{Corollary}
\newtheorem{remark}[thm]{Remark}
\newenvironment{proo}{{\bf Proof.}}{{\hfill $\square$}}

\title{On the moduli space of the Schwarzenberger bundles}
\author{Paolo Cascini}
\address{Dipartimento di Matematica\\ Viale Morgagni 67 A\\ 50134 Firenze\\ Italy}
\email{cascini@math.unifi.it}

\begin{document}
%\baselineskip= 1.7\baselineskip %(Pacific Version)
\begin{abstract}
For any odd $n$, we prove that the coherent sheaf $\mathcal F_A$
on $\PP^n_{\CC}$, defined as the cokernel of an injective map
$f:\OPN^{\oplus 2}\to \OPN(1)^{\oplus (n+2)}$, is
Mumford-Takemoto stable if and only if the map $f$ is stable, when
considered as a point of the projective space
$\PP(\Hom(\OPN(-1)^{\otimes 2},\OPN^{\otimes (n+2)})^*)$ under the
action  of the reductive group $\SL(2)\times\SL(n+2)$. This proves
a particular case of a conjecture of J.-M.Drezet and it implies
that  a component of the Maruyama scheme of the semi-stable
sheaves on $\PP^n$ of  rank $n$ and Chern polynomial $(1+t)^{n+2}$
is isomorphic to the Kronecher moduli $N(n+1,2,n+2)$, for any odd
$n$. In particular, such scheme defines a smooth minimal
compactification of the  moduli space of the rational normal
curves in $\PP^n$, that generalizes the construction defined by G.
Ellinsgrud, R. Piene and S. Str{\o}mme in the case $n=3$.
\end{abstract}

\maketitle

% \tableofcontents (indice)
%\vskip 1 cm

%%%%%%%%%%%%%%%%%%%%%%%%%%%%%%%%%%%%%%%%%%%%%%%%%%%%%%
\section{Introduction}
%%%%%%%%%%%%%%%%%%%%%%%%%%%%%%%%%%%%%%%%%%%%%%%%%%%%%%

Let us consider all the exact sequences:
\begin{equation}\label{suc.esatta}
0\longrightarrow I^*\otimes\mathcal O_{\PV}\stackrel {f_A}\longrightarrow
W^*\otimes\mathcal O_{\PV}(1)\longrightarrow \mathcal F_A\longrightarrow 0
\end{equation}
where $W$, $V$ and $I$ are complex vector spaces of
dimension $m+k$, $n+1$ and $k$ respectively,
$f_A$ is an injective morphism of sheaves canonically induced
by a linear map $A\in\PP(\Hom(W,I\otimes V)^*)$ 
($=\PP(\Hom(I^*\otimes \mathcal O_{\PV},W^*\otimes \mathcal O_{\PV}(1))^*)$)
and $\mathcal F_A =\coker f_A$ is a coherent sheaf of rank $m$ over
the projective space $\PV$ ($=(V^*\setminus\{0\})/\CC^*$).

In particular,  if $n=m$ and if the degeneracy locus of $f_A$ is empty,
then $\mathcal F_A$ is a vector bundle of rank $n$ on $\PP^n$, called
Steiner bundle. In \cite{GKZ} it is shown that $A$, considered as a
multidimensional matrix of size $(n+k)\times k\times (n+1)$, defines a
Steiner bundle $\mathcal F_A$ if and only if its hyperdeterminant does
 not vanish.

In \cite{AO}, the authors give a complete description of the
moduli space $S_{n,k}$ of the Steiner bundles on $\PP^n$: such moduli
space can be considered as an open subset of the categorical quotient:
$$\mathcal M_{n,m,k} = \PP(\Hom(W,I\otimes V)^*)//(\SL(I)\times\SL(W)).$$
with $n=m$. It is known that
$\mathcal M_{n,m,k}$ is canonically isomorphic to the Kronecker module
$N(n+1,k,m+k)$ defined as the quotient $\GG(W,I\otimes V)//\SL(I)$: the
isomorphism is given by considering the image $T_A:=A(W)$ of the
linear map $A:W\to I\otimes V$. Such
modules are extensively described in \cite{Dr1} and \cite{Dr2}. In
particular we have:
\begin{thm}\label{GIT.stabilita}
Let $A\in \PP(Hom(W,I\otimes V)^*)$ and $T=A(W)\subseteq I\otimes V$. The following are equivalent:
\begin{enumerate}
\item $A$ is semi-stable (resp. stable) under the action of
$\SL(I)\times\SL(W)$;
\item $T\in \GG(m+k,I\otimes V)$ is semi-stable (resp. stable) under
the action of $\SL(I)$;
\item for any non-empty subspace $I'\subsetneq I$
$$\frac {\dim T'} {\dim I'} \le \frac {\dim T} {\dim I}
\quad (\text{resp. }<)$$
where $T'=(I'\otimes V)\cap T$.
\end{enumerate}
\end{thm}

In general, if $m\ge n$, every element $A$ of $\mathcal M_{n,m,k}$ determines
a coherent sheaf $\mathcal F_A$ on $\PP^n=\PV$ of rank $m$: in fact,
every $A:W \rightarrow I\otimes V$ induces a morphism $f_A: I^*\otimes
\OPN\rightarrow W^*\otimes \OPN(1)$, as in (\ref{suc.esatta}).
We will call Steiner bundle of rank $m$, a vector bundle $\mathcal
F_A$ contained in the sequence (\ref{suc.esatta}) even when
$m\ge n$. Such bundles defines a moduli space $S_{n,m,k}$, that is an open
subset of the projective variety $\mathcal M_{n,m,k}$.

Important examples of rank $n$ Steiner bundles are the Schwarzenberger
bundles \cite{Sch}, defined by the morphism
$$f_A = \begin{pmatrix} %&x_0 &x_1 &\dots &x_n \\
                          x_0 &x_1   &\dots &x_n \\
                              &\ddots    &\ddots    & &\ddots      \\
                              &   &x_0  & x_1 & \dots   &x_n
        \end{pmatrix}^t
\in\PP(\Hom(I^*\otimes \OPN,W^*\otimes \OPN(1))^*).$$

The set of equivalence classes of these bundles is in one-one
correspondence with the variety $S_n$ of the rational normal curves.
In fact if $W(S)=\{H\in(\PP^n)^*|h^0((\mathcal F_A^*\otimes
\OPN(1))_{|H})\neq 0\}$ is
the scheme of the unstable hyperplanes of a Steiner bundle of rank $n$
$\mathcal F_A$, then $\mathcal F_A$ is a Schwarzenberger bundle if and
only if $W(S)$ is a rational normal curve in $(\PP^n)^*$ (see \cite{V}).

In particular, if $k=2$, all the indecomposable Steiner bundles are
Schwarzenberger
bundles (see \cite{DK}), and thus $S_n\simeq S_{n,2}\simeq \PP\GL(n+1)/\SL(2)$.
In this paper we will consider exactly this case. In fact we will
show that, if $k=2$ and $m$ is odd, then $A\in\PP(\Hom(W,I\otimes V)^*)$
 is stable if and only if the correspondent coherent sheaf $\mathcal F_A$ is
$\mu-$stable. This will imply the following:
\begin{thm}\label{main.thm}
$\mathcal M_{n,m,2}$ is isomorphic to the connected component of the Maru\-yama
moduli space  $\mathcal
M_{\PP^n}(m,c_1,\dots,c_n)$ containing the Steiner bundles. Such
component is smooth and irreducible.
\end{thm}

This result gives an affermative answer to a particular case  of a
question queried by J.-M. Drezet \cite{Dr3}. Before that,
R.M. Miro-Roig and G. Trautmann had
proved a similar result in the case $n=3$, $k=2$ and $m=3$
\cite{MT}.

Moreover the variety $\mathcal M_{n,n,2}$ defines a smooth
compactification of the moduli space of the rational normal curves in
$\PP^n$ for any odd $n$ (this result is proved in \cite{Dr2} and in \cite{ES}). In fact,
such construction generalizes  the one given in \cite{EPS}, defined as
the variety of nets of quadrics defining twisted cubics.

From a topological point of view, \cite{Dr2} provides a method to
compute the Betti numbers of
$\mathcal M_{n,m,2}$ (see also \cite{C} for further details).

\vskip 8 mm

I would like to thank V.Ancona and G.Ottaviani for many fruitful discussions 
and the referee for his very helpful comments.

\vskip 2 cm

%%%%%%%%%%%%%%%%%%%%%%%%%%%%%%%%%%%%%%%%%%%%%%%%%%%%%%
\section{Preliminares}
%%%%%%%%%%%%%%%%%%%%%%%%%%%%%%%%%%%%%%%%%%%%%%%%%%%%%%

Let $W, V$ and $I$ be complex vector spaces of dimension $m+2$, $n+1$
and $2$ respectively, with $2+m\le 2(n+1)$ and let us define $X=\PP(\Hom(W,I\otimes V)^*)$.

For any $\omega \in I$ we define $R_\omega=\omega\otimes V\subseteq
I\otimes V$: by theorem \ref{GIT.stabilita} we
have that an injective linear map $A\in X$ is
semi-stable (resp. stable) under the action of $\SL(I)\times\SL(W)$
if and only if, for any $\omega\in I$,
$$\dim R_\omega\cap T_A\le \frac{m+2} 2 \quad (\text{resp.} <),$$
where, we remind, $T_A$ is the image of $W$ by $A$ (the arithmetic assumption over $n$ and $m$ 
guarantees that $X^{ss}$ is not empty).

Let $D(A)$ denote the degeneracy locus of $f_A$, i.e. the set of all the points $x\in\PP^n$ such that
$\rank ((f_A)_x:I^*\otimes \mathcal O_{\PP^n,x}\to W^*\otimes\mathcal O_{\PP^n,x} (1)) \le 1$, then for
any $j\in \NN$ we construct the subsets:
\begin{eqnarray*}
S^j &=&\left\{A\in X^{ss}|\exists ~\omega\in I \text{ such that } \dim
R_\omega\cap T_A\ge j+m-n \right\}\qquad\text{and}\\
\tilde S^j&=&\{A\in X^{ss}|\dim D(A) \ge j - 2\}.
\end{eqnarray*}
These subsets canonically define two filtrations of $X$:
\begin{eqnarray*}
\emptyset= S^{j_0+1}\subseteq & S^{j_0}\subseteq \dots \subseteq
&S^2\subseteq S^1=X^{ss}\\
\emptyset\subseteq\dots\subseteq \tilde S^{j_0+1}\subseteq &\tilde
 S^{j_0}\subseteq \dots\subseteq &\tilde S^2\subseteq \tilde S^1=X^{ss}
\end{eqnarray*}
where $j_0=[\frac {m+3} 2]+n-m$ ($[x]$ denotes the integer part of
$x\in\QQ$).

It results $S^{j_0}=X^{ss}\setminus X^s$ and in particular it is empty if $m$ is odd.
Furthermore we have:
\begin{thm}\label{filtrazione}
\
\begin{enumerate}
\item $S^j\subseteq\tilde S^j\subseteq S^{j-1}$ for any $j\ge 2$;
\item $S^{2}=\tilde S^{2}$;
\item $S^{1}=\tilde S^{1}=X^{ss}$.
\end{enumerate}
In particular such subsets define a unique $G-$invariant filtration:
\begin{eqnarray*}
\emptyset=S^{j_0+1}\subseteq \tilde S^{j_0+1}\subseteq
S^{j_0}\subseteq \tilde S^{j_0}\subseteq\dots\\
\dots\subseteq S^3\subseteq
\tilde S^3\subseteq S^2 =\tilde S^2\subseteq S^1=\tilde S^1 =
X^{ss}
\end{eqnarray*}
\end{thm}

Before proving the theorem, we remind the following known lemma:
\begin{lemma}\label{thm2.8}
Let $F$ be a vector bundle of rank $f$ on a smooth projective variety
$X$ such that $c_{f-k+1}(F)\neq 0$ and let $\phi:\mathcal
O_X^k\longrightarrow F$ be a morphism with $k\le f$. Then the
degeneracy locus $D(\phi)=\{x\in X|\rank (\phi_x)\le k-1\}$ is
nonempty and $\codim D(\phi)\le f-k+1$.
\end{lemma}

\vskip 5 mm

\begin{proo} {\bf of theorem \ref{filtrazione}:}
\
\begin{enumerate}
\item Let $A\in S^j$, then there exists $\omega\in I$ such that
$\dim R_\omega\cap T_A\ge j + m - n$ and thus $\omega$ defines a
morphism of sheaves:
$\tilde {f_A}:\OPN\rightarrow\OPN(1)^{n - j + 2}$.
The degeneracy locus of  $\tilde f_A$ is contained in $D(A)$ and by
lemma \ref{thm2.8}, since
$c_{n -j +2}(\OPN(1)^{n - j +2})\neq 0$ if $j\ge 2$, it follows that
$\dim D(A)\ge j-2$, i.e.  $S^j\subseteq \tilde S^j$ for any $j\ge 2$.

Let now $A\in \tilde S^j$ and let us denote by $D_0(A)$ the variety of
all the points $x\in\PP^n$ such that $\rank (f_A)_x=0$.
We consider first the case
$D_0(A)\subsetneq D(A)$: each point 
$x\in\PP^n$ naturally defines an evaluation map $\eta_x:I\otimes V\to \PP(I^*)$.
Thus we can define $\pi:D(A)\setminus D_0(A)\rightarrow \PP(I^*)$ 
where $\pi(x)$ is the  only point of  $\eta_x(T_A)$ and,
since $\dim D(A)\setminus D_0(A)\ge j -2$, there exists $\omega\in I$
such that $\dim \pi^{-1}([\omega])\ge j -3$.

Let $R_\omega'=\{f\in I\otimes V|\eta_x(f) = [\omega] \text{ in } \PP(I^*) 
\text { for any }x\in\pi^{-1}([\omega])\}$: in order to compute
the dimension of $R_\omega'$
we consider $p_1,\dots,p_{j-2}\in\pi^{-1}([\omega])$ not contained in a
linear subspace $\PP^{j-4}\subseteq \PP^n$: such points define a
linear system of
$j - 2$ linearly independent equations whose solutions are contained
in $R'_\omega$ and thus we have
$\dim R'_\omega \le 2(n+1) - (j-2)= 2n + 4 -j$.

Since $T_A,R_\omega\subseteq R'_\omega$, we have that $$\dim T_A\cap
R_\omega \ge \dim T + \dim R_\omega - \dim R_\omega'\ge j + m - n - 1,$$
i.e. $A\in S^{j-1}$.

If $D_0(A)= D(A)$, then it can be similarly proven that for any
$\omega\in I$, $\dim R'_\omega\le 2n + 3 -j $ and thus $A\in
S^j\subseteq S^{j-1}$.

\vskip 4 mm

\item We have already proven that $S^2\subseteq \tilde S^2$. Let now $A\in
\tilde S^2$. As before, we can suppose $D_0(A)\subsetneq D(A)$.

Let $x\in D(A)\setminus D_0(A)$ and $\omega\in I$ such that
$\eta_x(T_A)=\{[\omega]\}$. If 
$R'_\omega=\{f\in I\otimes V|\eta_x(f) = [\omega] \text{ in } \PP(I^*)\}$,
then $\dim R'_\omega=2n +1$: $T_A,
R_\omega\subseteq R'_\omega$ and thus $\dim T_A\cap R_\omega \ge
(m+2)+(n +1)-(2n+1) = m -n +2 $, i.e. $A\in S^2$.

\vskip 4 mm

\item Both the equalities are trivial.
\end{enumerate}
\end{proo}

\vskip .7 cm

\begin{remark} In general $S^i\neq\tilde S^i$: let us consider, for
  instance, $n=m=3$ and
$$f_A=\begin{pmatrix} 0 &0& x_0& x_1& x_2 \cr
    x_0& x_1& 0& 0& x_3\end{pmatrix}^t.$$
Since $D(A)=\{(0:0:t_1:t_2)\}\simeq \PP^1$, $A\in\tilde S^3$; but
$S^3=\emptyset$ (see also prop. \ref{codimensione}).
\end{remark}

\vskip 1 cm

\begin{corol}\label{cor.codim}
If $m$ is odd and $A\in X^s=X^{ss}$ then $\codim D(A)\ge \frac {m+1}
2$.

If $m$ is even and $A\in X^{ss}$ (resp. $X^{s}$) then $\codim D(A)\ge
\frac m 2$ (resp. $>$).
\end{corol}

\begin{proo}
It suffices to notice that the previous theorem implies that $\tilde
S^{j_0+1}=\emptyset$ and that $S^{j_0}$ is the set of the properly
semi-stable points of $X$.
\end{proo}

\vskip .5 cm

\begin{prop}\label{codimensione}
If $m$ is odd, $A\in X$ is stable and $\codim D(A) = {\frac {m+1} 2}$, 
where $t={\frac {m+1} 2}$, 
  then, up to the action of
$\SL(I)\times\SL(W)\times\SL(V)$,
we have
$$f_A =\begin{pmatrix}x_0 &\dots & x_{t-1} & 0 &\dots &0 &x_t\\ 0
  &\dots &0 &x_0 &\dots &x_{t-1} &x_{t+1}
\end{pmatrix}^t.$$
\end{prop}
\begin{proof}
By the proof of theorem \ref{filtrazione} we have that for any
$\omega\in I$, $\dim(\omega\otimes V)\cap T_A \ge t$,
where, as before, $T_A$ is the image of $A$ as a subspace of $I\otimes V$, 
and in fact, by theorem \ref{GIT.stabilita} and since $A$ is stable, 
it results $\dim(\omega\otimes V)\cap T_A = t$. 

Thus we have, up to a change of basis,
$$f_A =\begin{pmatrix}f_0 &\dots & f_{t-1} & 0 &\dots &0 &f_t\\ 0
  &\dots &0 &g_0 &\dots &g_{t-1} &g_{t}
\end{pmatrix}^t,$$
where $<f_0,\dots, f_t>$ and $<g_0,\dots, g_t>$ are 
subspaces of $V$ of dimension $t+1$.

It is easily checked that $D(A) = V(f_0,\dots, f_t)\cup V(g_0, \dots,
g_t)\cup V(f_0,\dots,f_{t-1},
g_0,\dots,g_{t-1})$ and since $\codim D(A) = t$, it must be $\codim V(f_0,\dots,f_{t-1},
g_0,\dots,g_{t-1})= t$: this implies that $<f_0,\dots, f_{t-1}>=<g_0,\dots,g_{t-1}>$
and therefore we can assume $g_i=f_i$ for any $i=0,\dots,t-1$.

Moreover $g_t \notin <f_0,\dots,f_t>$  otherwise, up to the action of $\SL(I)\times\SL(W)$, 
it would be 
$$f_A =\begin{pmatrix}f_0 &\dots & f_{t-1} & 0 &\dots &0 &f_t\\ 0
  &\dots &0 &f_0 &\dots &f_{t-1} &0
\end{pmatrix}^t,$$
and by theorem \ref{GIT.stabilita}, $A$ would not be stable, because there 
would exist a vector $\omega\in I$ such that $\dim(\omega\otimes V)\cap T_A = t+1$.
Therefore $f_0,\dots,f_t$ are linearly independent and we can suppose
$f_i=x_i$ for some basis $\{x_0,\dots,x_n\}$ of $V$.  
\end{proof}

\vskip 2 cm

%%%%%%%%%%%%%%%%%%%%%%%%%%%%%%%%%%%%%%%%%%%%%%%%%%%%%%%%%%%%
\section{Proof of theorem \ref{main.thm}}
%%%%%%%%%%%%%%%%%%%%%%%%%%%%%%%%%%%%%%%%%%%%%%%%%%%%%%%%%%%%

For any coherent sheaf $\mathcal E$ of rank $r$ on $\PP^n$,
$\mathcal E_N$ will denote the normalized sheaf of $\mathcal E$,
i.e. $\mathcal E_N=\mathcal E(t_0)$ where $t_0\in\ZZ$ is such that
$-r< c_1(\mathcal E(t_0))\le 0$. Moreover $\hd(\mathcal E)$ will
be the homological dimension of $\mathcal E$ (cf. \cite{OSS}) and
$S(\mathcal E)$ the singular locus of $\mathcal E$, i.e.
$S(\mathcal E)=\{x|\dim \mathcal E_x>r\}$.

\vskip 5 mm

In this section we will only consider sheaves $\mathcal F_A$ of
odd rank $m$, i.e. such that $(c_1(\mathcal F_A), m)=1$. Hence the
Mumford-Takemoto stability (also said $\mu-$stability) of these
sheaves coincides with their Gieseker stability. Moreover
$\mathcal F_A$ is stable if and only if it is semi-stable. Thus,
before proceeding with the proof of theorem \ref{main.thm}, we are
interested to study the relation between the G.I.T. stability of
maps and the $\mu-$stability  of their cokernels.

In fact we have:
\begin{thm}\label{stabilita}
Let $k=2$ and $m\in\NN$ odd. Then the following are equivalent:
\begin{enumerate}
\item $T_A\in \GG(m+2,I\otimes V)$ is G.I.T. stable;
\item $\mathcal F_A$ is $\mu$-stable.
\end{enumerate}
\end{thm}
The main tool needed to prove the theorem is the following lemma:
\begin{lemma}\label{main.lemma}
Let $A\in \PP(\Hom(W,I\otimes V)^*)$ be a stable map, then
\begin{equation}\label{hoppe}
\HH^0((\wedge^r \mathcal F_A)^{**}_N)=0
\end{equation}
for any $r=1,\dots,m-1$.
\end{lemma}
Later on, we will show that the vanishing of the cohomology groups
in (\ref{hoppe}) will imply the $\mu$-stability of the sheaf
$\mathcal F_A$.

\vskip .8 cm

Before proceeding with the proof of lemma \ref{main.lemma},
we want to recall some facts that will be useful during the proof:
although many of these results are well known, we report them for
completeness.

For the proof of the following two propositions, see \cite{HL}
prop. 1.1.6 and prop. 1.1.10:

\begin{prop}\label{reflexive1}
Let $E$ be a coherent sheaf of codimension $c$ on a smooth projective
variety $Z$. Then the sheaves $\Extl^q(E,\omega_Z)$ are supported on
$\Supp(E)$ and $\Extl^q(E,\omega_Z)=0$ for all $q<c$.
\end{prop}

\begin{prop}\label{reflexive}
Let $E$ be a coherent sheaf on a smooth projective variety $Z$.
Then the following conditions are equivalent:
\begin{enumerate}
\item $\codim(\Extl^q(E,\omega_Z))\ge q+1$ for any $q\ge 1$;
\item the canonical map $E \rightarrow E^{**}$ is injective.
\end{enumerate}
Similarly, the following are equivalent:
\begin{enumerate}
\item $\codim(\Extl^q(E,\omega_Z))\ge q+2$ for any $q\ge 1$;
\item $E$ is the dual of a coherent sheaf;
\item $E$ is reflexive.
\end{enumerate}
\end{prop}

We will also need:

\begin{lemma}\label{koszul}
Let $s$ be a section of a vector bundle $E$ of rank
$r$ on an algebraic variety $Z$ and let $Z_0$ be the zero locus of
$s$. If $Z_0$ is of codimension $r'\le r$
%and $Z_0\subseteq
%Z'\subseteq Z$, where $Z'$ is a complete intersection subvariety of
%$Z$ of codimension $r'$,
then the Koszul complex associated to $s$ induces an exact
sequence of the first $r'+1$ terms:
\begin{equation}\label{koszul.complex}
0\rightarrow \det E^*\rightarrow\wedge^{r-1} E^*\rightarrow
\dots\rightarrow \wedge^{r-r'}E^*
\end{equation}
\end{lemma}
\begin{proo}
By Bertini theorem, it easily follows that there exists a complete
intersection subvariety $Z'\subseteq Z$ of  codimension $r'$ and containing $Z_0$.

It is enough to prove the lemma after restricting the bundle $E^*$
into a trivializing open subset $U\subseteq Z$ such that
$E^*_U\simeq \CC^r \otimes \mathcal O_U$ where $\CC^r$ is spanned
by $e_1,\dots, e_r$. We can suppose that, with respect of this
frame, $s=(s_1,\dots,s_r)$ and that $Z'\cap U$ is the zero locus
of $s'=(s_1,\dots,s_{r'})$.

Let us proceed by induction on $r-r'$.
If $r=r'$, then $Z_0$ is a complete intersection and the Koszul complex
 is exact (\cite{GH}, pag.688).

Let us suppose now $r>r'$ and let $E_k=\wedge^k E_U^*$ and
$F_k=\wedge^k \CC^{r-1}\otimes \mathcal O_U\subseteq
\wedge^k \CC^r\otimes \mathcal O_U\simeq E_k$.
The quotient $Q_k$ of $E_k$ by $F_k$ is isomorphic to
$(e_r\otimes\wedge^{k-1}\CC^{r-1})\otimes \mathcal O_U$, moreover the map
$\delta:E_k\to E_{k-1}$ induces, in a canonical way, two maps
$\delta':F_k\to F_{k-1}$ and $\delta'':Q_k\to Q_{k-1}$ so that $F_*$
and $Q_*$ are again Koszul complexes contained in the
commutative diagram:

\begin{equation*}
\begin{CD}
      & &  &         &0            &&0             &&&&0              \\
      & &  &         &\downarrow   &&\downarrow    &&&&\downarrow     \\
      &  &0   @>>>    F_{r-1} @>>>  F_{r-2} @>>> \dots @>>> F_{r-r'}  \\
      & &\downarrow  &&\downarrow  &&\downarrow    &&&&\downarrow     \\
0 @>>> E_r @>>>    E_{r-1}  @>>>  E_{r-2} @>>> \dots @>>> E_{r-r'}    \\
      & &\downarrow  &&\downarrow  &&\downarrow    &&&&\downarrow     \\
0 @>>> Q_r @>>>    Q_{r-1}  @>>>  Q_{r-2} @>>> \dots @>>> Q_{r-r'}    \\
       &     &\downarrow  &&\downarrow &&\downarrow  &&&&\downarrow   \\
       &            &0            &&0            &&0 &&&& 0
\end{CD}
\end{equation*}

\vskip 5 mm

By induction hypothesis, $\HH_*(F_*)=\HH_*(Q_*)=0$. Thus also
$\HH_*(E_*)=0$, i.e. the sequence (\ref{koszul.complex}) is exact.
\end{proo}

\vskip 5 mm

\begin{lemma}\label{successione.lunga}
For any $r=1, \dots, m-1$, the sheaf $\wedge^r \mathcal F_A$ is
contained in the exact sequence:
\begin{equation}\label{sel.eis}
  I^*\otimes\wedge^{r-1}W^*\otimes\OPN(r-1)\rightarrow\wedge^r
  W^*\otimes\OPN(r)\rightarrow \wedge^r \mathcal F_A\rightarrow 0.
\end{equation}
Moreover, if $A\in X$ is G.I.T. stable and $r\le \frac{m-1} 2$ then
the sequence
\begin{multline}\label{sel}
  0\rightarrow S^r I^*\otimes \OPN \rightarrow S^{r-1}I^*\otimes W^*\otimes
  \OPN(1)\rightarrow\dots\\
  \dots\rightarrow
  I^*\otimes\wedge^{r-1}W^*\otimes\OPN(r-1)\rightarrow\wedge^r
  W^*\otimes\OPN(r)\rightarrow \wedge^r \mathcal F_A\rightarrow 0
\end{multline}
is exact. In particular, $\hd(\wedge^r
 \mathcal F_A) \le r$.
\end{lemma}

\begin{proof}
The exactness of (\ref{sel.eis}) is proven in
(\cite{Eis}, pag. 571).

In order to prove the exactness of (\ref{sel}), we proceed by mimicking
the proof of the existence of the Eagon-Northcott complex given in \cite{GP}.

Let $Z=\PV\times\PP(I^*)$ and let $\pi:Z\rightarrow \PV$ be the
projection onto the first space.
The morphism $A$ defines a section $a:{\mathcal O}_Z\rightarrow
W^*\otimes {\mathcal O}_Z(1,1)$, given by $a=(y_0 f_{0,i} +y_1 f_{1,i}
)_{i=1}^{m+2}$ where the $f_{i,j}$'s are the entries of $A$
and $y_0,y_1$ are the coordinates of $\PP(I^*)$.

The zero locus of $a$ is
$\tilde Z = \cap_{i} V(y_0 f_{0,i} +y_1 f_{1,i})\subseteq Z,$
and the Koszul complex associated is given by:
\begin{eqnarray*}
0\rightarrow\wedge^{m+2}W\otimes {\mathcal
  O}_Z(-m-2,-m-2)\rightarrow\dots&\\
 \dots \rightarrow
  \wedge^2W\otimes {\mathcal O}_Z(-2,-2)\rightarrow
  W\otimes  {\mathcal O}_Z(-1,-1)&\rightarrow
  {\mathcal O}_Z\rightarrow {\mathcal O}_{\tilde Z}\rightarrow 0
\end{eqnarray*}
We have that  $\pi(\tilde Z)\subseteq D(A)$ and, since $A$ is
stable, then $\dim \tilde Z\le \dim D(A)+1\le n-\frac {m-1}
2$ (Cor. \ref{cor.codim}).

By lemma \ref{koszul}, the sequence:
$$0\rightarrow\wedge^{m+2}W\otimes {\mathcal
  O}_Z(-r-2,-r-2)\rightarrow \dots \rightarrow
 \wedge^{m+2-r}W\otimes {\mathcal O}_Z(-2,-2).$$
is exact for any $r\le \frac{m-1} 2$.

Since each fiber of $\pi$ is isomorphic to $\PP(I^*)$,
it results:
$$R^i\pi_*(\mathcal O_Z(-2-j,-2-j))=\begin{cases}
  S^jI^*\otimes\OPN(-2-j)\qquad&\text{ if } i=1\quad j\ge 0 \\
0 \qquad &\text{ if } i\neq 1\quad j\ge 0 \end{cases} $$
where $R^i\pi_*$ is the higher direct image functor associated to
$\pi$ (see \cite{Har}, ch. III, 8). Moreover
 $\wedge^{m+2-r+j} W\simeq \wedge^{r-j} W^*$,
that yields the exact sequence:
\begin{equation}\label{sel1}
0\rightarrow S^rI^* \otimes
\OPN(-r-2)\rightarrow \dots \rightarrow \wedge^r W^*\otimes\OPN(-2).
\end{equation}
The exactness of the sequence (\ref{sel}) follows by gluing
(\ref{sel1}) tensored by $\OPN(r+2)$ with (\ref{sel.eis}):
in fact, in both the sequences the morphisms
$$I^*\otimes\wedge^{r-1}W^*\otimes\OPN(r-1)\rightarrow\wedge^r
W^*\otimes\OPN(r)$$
are canonically defined.
\end{proof}

\vskip 4 mm
\begin{lemma}\label{riflessivita}
Let $E$ be a coherent sheaf and $S(E)$ its singular locus. If $\codim
S(E)\ge \hd(E)+2$, then $E$ is  reflexive.
\end{lemma}
\begin{proof}
Let $t=\hd(E)$ and let us consider a resolution of $E$:
$$0\rightarrow F_t\rightarrow F_{t-1}\rightarrow F_{t-2}\rightarrow\dots
\rightarrow F_2\rightarrow F_1\rightarrow F_0 \rightarrow E\rightarrow 0$$
where each $F_i$ is a direct sum of line bundles.

Let us split the sequence into short exact sequences:
$$0\rightarrow F_t \rightarrow F_{t-1}\rightarrow
G_{t-1}\rightarrow 0$$ $$0\rightarrow G_i \rightarrow
F_{i-1}\rightarrow G_{i-1}\rightarrow 0$$ $$0\rightarrow G_1
\rightarrow F_{0}\rightarrow E\rightarrow 0$$ then, applying the
functor $\Extl^i(\cdot,\omega_{\PP^n})$, we get:
$$\Extl^i(E,\omega_{\PP^n})=\Extl^{i-1}(G_1,\omega_{\PP^n})=
%\Extl^{i-2}(G_2,\omega_{\PP^n})=
\dots=\Extl^{i-t+1}(G_{t-1},\omega_{\PP^n})=0,$$
for any $i$ such that $i-t+1\ge 2$.
Thus $\Extl^i(E,\omega_{\PP^n})=0$ for any $i\ge t+1$.

Moreover $$\codim \Extl^i(E,\omega_{\PP^n})\ge \codim S(E)\ge t+2\ge i+2$$
if $1\le i\le t$ and thus by prop. \ref{reflexive}, $E$ is reflexive.
\end{proof}

\vskip 4 mm

\begin{lemma}\label{dualita}
Let $E$ be a sheaf of rank $m$ on $\PP^n$ such that $\codim S(E)\ge
2$. Then, for any $r=1,\dots, m-1$, we have
\begin{equation}\label{ciao}
(\wedge^r E)^{**} = (\wedge^{m-r} E)^*\otimes
\OPN(c_1(E))
\end{equation}
\end{lemma}

\begin{proof}
The injective map $\wedge^{m-r}E\to \Hom(\wedge^{r}E,\wedge^m E)$ induces
the exact sequence:
$$0\rightarrow\wedge^{m-r}E
   \rightarrow\Hom(\wedge^{r}E,\wedge^m E)
   \rightarrow\mathcal E\rightarrow 0,$$ 
where $\mathcal E$ is a $0$-rank sheaf such that $\codim\mathcal E\ge\codim S(E)\ge 2$:
thus, dualizing this sequence and observing that, by prop.
\ref{reflexive1}, $\mathcal E^*=\Extl^1(\mathcal
E,\OPN)=0$, it results $(\wedge^{m-r}E)^*\simeq \Hom(\wedge^{r}E,\wedge^m E)^*
\simeq(\wedge^r E)^{**}\otimes (\wedge^m E)^*$, in fact by proposition 
\ref{reflexive} all these sheaves are torsion-free. Moreover $(\wedge^m E)^{**}\simeq 
\OPN(c_1(E))$ and therefore (\ref{ciao}) follows.
\end{proof}

\vskip 1.5 cm

Thus it results $\HH^0((\wedge^r \mathcal F_A)_N^{**})=
\HH^0((\wedge^{m-r}\mathcal F_A)^*(t_0))$ for suitable $t_0\in
\ZZ$: we want to prove that such cohomology group is null.

Let us distinguish 2 cases:

\vskip 4 mm

I. $r\ge {\frac {m+1} 2}$:

Let $t_0$ be such that $(\wedge^r \mathcal F_A)_N^{**} =
(\wedge^{m-r} \mathcal F_A)^*(t_0)$. By the sequence
(\ref{sel.eis}), we have: 
$$0\rightarrow (\wedge^{m-r}\mathcal
F_A)^*(t_0)\rightarrow \wedge^{m-r} W\otimes \OPN \stackrel B
\longrightarrow \wedge^{m-r-1}W\otimes I^\otimes \OPN(1).$$
Thus if $\HH^0((\wedge^{m-r}\mathcal F_A)^*(t_0))\neq 0$ then
there exists $b:\OPN\hookrightarrow \wedge^{m-r} W\otimes \OPN$
(i.e. $b\in \wedge^{m-r} W$) such that $B\circ b=0$.

It is easy to see that if $A$ is injective then $B\circ b=0$ implies $b=0$.

\vskip 8 mm
II. $r\le {\frac {m-1} 2}$:

If4 $\codim D(A)\ge 2+{\frac {m-1} 2}={\frac {m+3} 2}$, then, by
lemma \ref{riflessivita}, the sheaf $\wedge^r\mathcal F$ is
reflexive, and by the sequence (\ref{sel}), it is easy to show
that: $$\HH^0((\wedge^r \mathcal F_A)^{**}_N)=\HH^0((\wedge^r
\mathcal F_A)_N)=0$$
  
By corollary \ref{cor.codim}, we have that the G.I.T. stability of
$A$ implies that $\codim D(A)\ge {\frac {m+1} 2}$ and thus it just
remains to consider the matrices $A$ such that $\codim D(A)=
{\frac {m+1} 2}$.

\begin{lemma} If $\codim D(A)=\frac {m+1} 2$
then $\HH^0((\wedge^r \mathcal F_A)^{**}_N)=0$ for any
$r=1,\dots,m-1$.
\end{lemma}
\begin{proof}
By prop. \ref{codimensione}, we can suppose
$$f_A =\begin{pmatrix}x_0 &\dots & x_{t-1} & 0 &\dots &0 &x_t\\ 0
  &\dots &0 &x_0 &\dots &x_{t-1} &x_{t+1}
\end{pmatrix}^t.$$
Moreover the same technique used above can be applied to prove the
thesis for all $r\neq t$.

Thus it suffices to show that $\HH^0((\wedge^t\mathcal
F_A)^*(t_0))= \HH^0((\wedge^{t-1}\mathcal F_A)^{**}_N)=0$. It is
easily checked that by dualizing the sequence (\ref{sel}) we get
the sequence $$0\rightarrow (\wedge^{t}\mathcal
F_A)^*(t_0)\rightarrow \wedge^t W\otimes \OPN(1) \stackrel B
\longrightarrow \wedge^{t-1}W \otimes I \otimes \OPN(2).$$
Thus we just need to prove that if $b:\OPN\to\wedge^t W\otimes
\OPN(1)$ is such that $B\circ b=0$ then $b=0$: this is a direct
computation.
\end{proof}

\vskip 5 mm

Thus lemma \ref{main.lemma} is completely proven. We can proceed
now with the proof of theorem \ref{stabilita}.

\vskip 5 mm
\begin{proof}[Proof of theorem \ref{stabilita}]
The first statement of prop. \ref{reflexive} easily implies that
the map $\mathcal F_A \rightarrow (\mathcal F_A)^{**}$ is
injective. Since $(\mathcal F_A)^{**}$ is torsion-free, so is
$\mathcal F_A$.

Let now $\mathcal E\subseteq \mathcal F_A$ be a torsion-free sub-sheaf
of rank $r$.
then $\OPN(c_1(\mathcal E))=(\wedge^r\mathcal E)^{**}\subseteq
(\wedge^r \mathcal F_A)^{**}$: since $\HH^0((\wedge^r \mathcal
F_A)^{**}_N)=0$ (lemma \ref{main.lemma}), it results $c_1(\mathcal E)<
\mu(\wedge^r \mathcal F_A)=r\mu(\mathcal F_A)$, i.e. $\mu(\mathcal
E)<\mu(\mathcal F_A)$.
Thus $\mathcal F_A$ is $\mu$-stable.

Vice-versa, let $A\in X$  be a non-stable
matrix. Then, by theorem \ref{GIT.stabilita}, we can write
$A=\begin{pmatrix}0 & A_0\\ A_1 &
                      A_2 \end{pmatrix}$
where $A_0$ is a vector of length $s\le\frac{m+1} 2$.
Thus $A_0$ defines a sub-sheaf $\mathcal F_{A_0}\subseteq \mathcal F_A$
that is  contained in the exact sequence:
$$0\longrightarrow\OPN\stackrel {f_{A_0}}\longrightarrow\OPN(1)^s
\longrightarrow \mathcal F_{A_0}\longrightarrow 0.$$
It is easily checked that $\mu(\mathcal F_{A_0})>\mu(\mathcal F_{A})$.
\end{proof}

\vskip 5 mm

We show now an interesting relation within the automorphism group of
$\mathcal F_A$ and the stabilizer of $A$.

\begin{thm}\label{simple}
 Let $A\in X$ such that $\mathcal F_A$ is simple, i.e.
  $\Aut(\mathcal F_A)=\CC^*$. Then
\begin{equation}\label{stabilizer}
\Stab_G(A)=\{(\lambda \Id_{2},\mu \Id_{m+2})\in
G|\lambda^{n+k}=\mu^k=1\}.
\end{equation}
In particular
  $\dim \Stab_G(A)=0$ and $\dim \mathcal M_{n,m,2}=\dim X - \dim G$, for any
  $m<2n$.
\end{thm}

\begin{proof}
Let us prove first that any $f\in \Aut(\mathcal F_A)$ is uniquely
determined by a morphism of sequences:

$$\begin{CD}
0@>>> I^*\otimes\OPN @>{f_A}>> W^*\otimes \OPN(1) @>>> \mathcal F_A @>>> 0\\
    &&   @VPVV            @VQVV               @VVfV                   \\
0@>>> I^*\otimes\OPN @>{f_A}>> W^*\otimes\OPN(1) @>>> \mathcal F_A @>>> 0 \\
\end{CD}$$

This is a direct consequence of the fact that, by the vanishing of
$\Hom(\OPN^{n+k},\OPN(-1)^k)$ and $\Ext^1(\OPN^{n+k},\OPN(-1)^k)$ we
get the exact sequence:
$$0\rightarrow\Hom(\OPN^{n+k},\OPN^{n+k})
\rightarrow\Hom(\OPN^{n+k},\mathcal F_A)\rightarrow 0.$$
Thus if $\mathcal F_A$ is simple, then the only automorphisms of
$\mathcal F_A$ are the homotheties, that implies (\ref{stabilizer}).
\end{proof}

\vskip .5 cm

We are ready now to prove theorem \ref{main.thm}. Let $c_1,\dots,
c_n$ be the Chern classes of $\mathcal F_A$ and let $\mathcal
M_{\PP^n}(m,c_1,\dots,c_n)$ be the Maruyama moduli space of all
the $\mu$-stable sheaves or rank $m$ and Chern classes $c_1,\dots,
c_n$. By theorem \ref{stabilita}, if $m$ is odd, each $A\in
\mathcal M_{n,m,2}$ defines uniquely an isomorphism class of
coherent sheaves $[\mathcal F_A]\in \mathcal
M_{\PP^n}(m,c_1,\dots,c_n)$ and thus there exists an injective
projective morphism: $$\phi:\mathcal M_{n,m,2}\longrightarrow
\mathcal M_{\PP^n}(m,c_1,\dots,c_n)$$ Moreover, by the sequence
(\ref{suc.esatta}) that defines $\mathcal F_A$, it is easily
checked that $\Ext^2(\mathcal F_A, \mathcal F_A)=0$, i.e. every
point of the image of $\phi$ is a smooth point of the Maruyama
moduli space. By theorem \ref{simple} and by sequence
(\ref{suc.esatta}), it immediately follows that $\dim \mathcal
M_{n,m,2}=\dim \Ext^1(\mathcal F_A,\mathcal F_A)=\dim_{[\mathcal
F_A]}\mathcal M_{\PP^n}(m,c_1,\dots,c_n)$, and thus, by Stein
factorization theorem, $\phi$ maps isomorphically $\mathcal
M_{n,m,2}$ onto a smooth connected component of $\mathcal
M_{\PP^n}(m,c_1,\dots,c_n)$. This completely proves theorem
\ref{main.thm}.

\vskip 5 mm

\begin{remark}
In general, the same result does not hold for higher $k$, even if $(m,k)=1$ (i.e. in the case where all the semi-stable sheaves are stable).
Consider, for istance, the matrix
$$f=\begin{pmatrix} x_1 &x_0
&0&0&0\\0&x_2&x_1&x_0&0\\0&0&0&x_2&x_1\end{pmatrix}^t $$
and the sequence
$$0\longrightarrow \mathcal O_{\PP^2}^3\stackrel f\longrightarrow
   \mathcal O_{\PP^2}(1)^5\longrightarrow
   \mathcal F_A\longrightarrow 0.$$
Then, since the degeneracy locus of the sheaf $\mathcal
F$ , $\{x_1=0\}$, is of codimension $1$, $\mathcal F$ is not
torsion-free and in
particular it cannot be $\mu$-stable. On the other hand, using
theorem \ref{GIT.stabilita}, it is a direct computation to prove that
$f$ defines a stable morphism.
\end{remark}

\vskip 2 cm

%%%%%%%%%%%%%%%%%%%%%%%%%%%%%%%%%%%%%%%%%%%%%%%%%%%%%%
\providecommand{\bysame}{\leavevmode\hbox to3em{\hrulefill}}

\end{document}